\documentstyle [12pt]{article}
\renewcommand{\H}{{\cal H}}

\newcommand{\M}{{\cal M}}
\newcommand{\K}{{\cal K}}
\newcommand{\B}{{\cal B}}

\newcommand{\A}{{\cal A}}
\newcommand{\F}{{\cal F}}
\newcommand{\U}{{\cal U}}

\begin {document}
\sloppy
\large

\begin {center}
{\bf ON COCYCLE CONJUGACY OF QUASIFREE}
\end {center}
\begin {center}
{\bf ENDOMORPHISMS SEMIGROUPS}
\end {center}
\begin {center}
{\bf ON THE CAR ALGEBRA}
\end {center}

\bigskip

\begin {center}
{\bf GRIGORI G. AMOSOV}
\end {center}

\begin {center}
{\bf Department of Higher Mathematics, Moscow Institute of Physics and
Technology, 141700~Dolgoprudni, Russia, E-mail: gramos@mail.sitek.ru}
\end {center}

   W. Arveson has described a cocycle conjugacy class $\U(\alpha )$ of
$E_0$-semigroup $\alpha $ on $\B(\H)$ which is
a factor of type $\rm I$. Under some
conditions on $\alpha $ there is a $E_0$-semigroup
$\beta \in \U(\alpha )$ being a flow of shifts in the sense of
R.T.Powers (see [11]). We study quasifree endomorphisms
semigroups $\alpha $ on the hyperfinite factor $\M=\pi (\A(\K))''$
generated by
the representation $\pi $ of the algebra of the canonical anticommutation
relations $\A(\K)$ over a separable Hilbert space $\K$.
The type of $\M$ can be $\rm I$, $\rm II$ or $\rm III$ depending on $\pi $.
The cocycle conjugacy class $\U(\alpha )$ is described in the terms of
initial isometrical semigroup in $\K$ and an analogue of the Arveson result
for the hyperfinite factor $\M$ of type $II_1$ and $III_{\lambda },
0<\lambda <1,$ is introduced.

\medskip

{\bf 1. Introduction.} Let $\cal M$ be the $W^*$-algebra acting in a Hilbert
space. One-parameter unital semigroup $\alpha _t\in End(\M),\
t\geq 0,$
is called a $E_0$-semigroup if
every function $\eta (\alpha _t(x))$ is continuous in $t$ for
$x\in \cal M$ and $\eta \in \cal M\it _*$.
Given a $E_0$-semigroup $\alpha $ one can define its generator
$\delta (a)=lim_{t\to 0}\frac {\alpha _t(a)-a}{t}$ for
$a\in dom\delta$, where $dom\delta$ is $\sigma $-week
dense in $\cal M$ (see [14]).
Two $E_0$-semigroups $\alpha $ and $\beta $ are called to be cocycle conjugate
if there is a strong continuous family of unitaries $U_t\in \cal M,\
\it t\geq 0,$ named a cocycle, such that
$\beta _t(\cdot )=U_t\alpha _t(\cdot )U_t^*,\ U_{t+s}=U_t\alpha _t(U_s),\
t,s\geq 0$ (see [11,14,15]). Notice that if semigroups $\alpha $
and $\beta $ have generators differ on a bounded derivation, then
$\alpha $ and $\beta $ are cocycle conjugate (see [15]).
This case is associated with the differentiable cocycle $(U_t)_{t\geq 0}$
and it is not the general one. The discussion on the cocycle conjugacy of
automorphisms semigroups on the
$W^*$-algebra $\cal M=B(H)$ one can see in monograph [14].
A notion of the cocycle conjugacy of endomorphisms semigroups
on the $W^*$-algebra $\cal M$ was given by W. Arveson.
He studied the case of
$\cal M=B(H)$ in [11].

     Let $\A=\A(\K)$ be the $C^*$-algebra of the canonical anticommutation
relations (CAR) over a Hilbert space $\K$. It means that there
is a map $f\to a(f)$ from $\K$ to the $C^*$-algebra $\cal A$ (with
the unit $\bf 1$, the involution * and the norm $||\cdot ||$) satisfying
the following properties:

     $1)\ a(\lambda f+g)=\lambda ^*a(f)+a(g)$ for all $f,g\in \K,\
\lambda \in \bf C$,

      2) (CAR) $a(f)a(g)+a(g)a(f)=0$,

$\ \ \ \ \ \ \ \ \ \ \ \ a^*(f)a(g)+a(g)a^*(f)=(g,f)\bf 1$,

      3) the polinomials in all $a(f),a^*(g)$ are dense in $\cal A$ by the norm,
$||a(f)||=||f||_{\K}$.

      Every state on $\cal A$ that is a positive linear functional
$\phi \in \cal A\it ^*,\ \phi ({\bf 1})=1$, is determined by its values on
(Wick) normal
ordered monomials $a^*(f_1)...a^*(f_m)a(g_1)...a(g_n)$. The operator
$R,\ 0<R<1,$ in $\cal B(\K)$ determines a state $\omega _R$ satisfying
the condition
$$
\omega _R(a^*(f_m)...a^*(f_1)a(g_1)...a(g_n))=\delta _{nm}det((f_i,Rg_j)).
$$
Such $\omega _R$ is called a quasifree state.
Let us define a representation $\pi _R$ of the $C^*$-algebra $\cal A$ in a Hilbert
space $\cal H=F(\K)\otimes \cal F(\K)$ by the formula
$$
\pi _R(a(f))=a((1-R)^{1/2}f)\otimes \rm \Gamma \it +1\otimes a^*(JR^{1/2}f),\
f\in {\cal K},
$$
$$
\pi _R(\bf 1)=\rm Id,\ \ \ \ \ \ \ \ \ \ \ \ \ \ \ \ \ \ \ \ \ \ \ \ \ \ \ \ \ \ \
\ \ \ \ \ \ \ \ \ \ \ \ \ \ \ \ \ \ \ \ \ \ \ \ \ \ \
$$
where $\cal F(\K)$ is the antisymmetric (fermion) Fock space over $K$ with
a vacuum vector $\Omega $, $J$ is some antiunitary in $K$, $\Gamma $ is
a hermitian unitary operator completely defined by the condition $\Gamma a(f)=
-a(f)\Gamma ,f\in h,\Gamma \Omega =\Omega $.
In this representation the state $\omega _R$ becomes the vector one,
$\omega _R(x)=(\Omega \otimes \Omega ,\pi _R(x)\Omega \otimes \Omega ),\
x\in \A$. Therefore the triple $(\pi _R,\cal H\it _R,\Omega \otimes \Omega )$,
where $\cal H\it _R=\overline {\pi _R(\A)\Omega \otimes \Omega }$,
is the Gelfand-Neumark-Segal (GNS) representation of the $C^*$-algebra $\cal A$
associated with the state $\omega _R$.
The state $\omega _R$ is pure if and only if $R=P,P^2=P$ is
an orthogonal projection.
In this case the GNS representation $\pi _P$
acting in a Hilbert space $\cal H=\cal F(\it (I-P)\K)\otimes \cal F(\it JPJ\K)$
yelds the $W^*$-algebra $\cal M\it _P(\K)=\pi _P(\cal A)\it ''=\cal
B(H)$. Setting a twopoints function of the state $\omega $ to be
$$
\omega (a^*(f)a(g))=\nu (f,g),\  f,g\in \K,
$$
and fixing a numerical parameter $\nu \in [0,1/2]$ one can consider
the quasifree state
$\omega =\omega _{\nu }$ assosiated with the operator $R=\nu \bf 1$ in $\K$.
If $\nu \neq 0$, then $\omega $ is exact ($\omega (x^*x)=0$
implies $x=0$). In the GNS representation $\pi _{\nu }$ acting
in a Hilbert space $\cal H=\cal F(\K)\otimes \cal F(\K)$
the $C^*$-algebra $\pi _{\nu }(\cal A)$ generates the $W^*$-algebra
$\cal M\it _{\nu }=\cal M\it _{\nu }(\K)=\pi _{\nu }(\cal A(\K))''$.
If $0<\nu <1/2$, then $\cal M\it _{\nu }$ is a hyperfinite
factor of type $III_{\lambda }$,
where $\lambda =\nu /(1-\nu )$. In the case of $\nu =1/2$ the state
$\omega $ is a trace and $\cal M\it _{\nu }$ is a hyperfinite factor of
type $II_1$. The case of $\nu =0$ is associated with the vacuum state and
the Fock representation with $\cal M\it _0=\cal B(F(\K))$ (see [9-10,16-20]).

   Let a quasifree endomorphism $\alpha $ on the algebra $\pi _{R}(\cal A$)
act on the generating elements by the formula
$\alpha (\pi _{R}(a(f)))=\pi _{R}(a(Vf)),\ f\in \K$,
where $V$ is an isometry in a Hilbert space $K$ commuting with $R$.
If $kerR=ker(I-R)=0$ the state $\omega _{R}$ is exact and the vector
$\Omega \otimes \Omega $ is separating for the $W^*$-factor
$\cal M\it _{R}=\pi _R(\A)''$.
It allowes to show that
$\alpha $ can be extended
to a quasifree endomorphism of $\cal M\it _{R}$ (see [16]). We
denote this endomorpism by $B_{R}(V)$.
The procedure of the passage from an operator in a Hilbert space $\K$ to a map
on the algebra ${\cal M}_{R}$ is called the (quasifree) lifting.
If $(V_t)_{t\geq 0}$ is a $C_0$-semigroup of isometries in $\K$
commuting with $R$, then
the semigroup $(B_{R}(V_t))_{t\geq 0}$ is a $E_0$-semigroup on
the $W^*$-factor ${\cal M}_{R}$.
We call these semigroups the quasifree ones (see [9-10,16-19]).

The endomorphism $\alpha$ of $W^*$-algebra $\cal M$ is called
a shift if $\cap _{n=1}^{+\infty }\alpha ^n(\cal M)=\bf C1$.
The $E_0$-semigroup $(\alpha _t)_{t\geq 0}$ is called a flow of shifts if
$\alpha _t$ is a shift for every fix $t>0$.
In [12] R.T. Powers introduced a flow of shifts on the $W^*$-algebra
$\cal M\it _0=\cal B(F(\K))$.
It was obtained by an extension of the semigroup
$(\alpha _t)_{t\geq 0}$ acting
on the generating elements of the $C^*$-algebra $\pi _0({\cal A}({\cal K})),\
{\cal K}=L_2(0,+\infty ),$ by the formula
$\alpha _t(\pi _0(a(f)))=\pi _0(a(S_tf)),\ t\geq 0,\ f\in \K$, where $(S_t)_{t\geq 0}$ is
a $C_0$-semigroup of right shifts in $\K$ defined by the formula
$(S_tf)(x)=f(x-t)$ for $x>t$ and $(S_tf)(x)=0$
for $0<x<t,\ f\in \cal K$. In [13] it is asserted that the quasifree
$E_0$-semigroup $(B_{\nu }(S_t))_{t\geq 0}$ on the hyperfinite factor
${\cal M}_{\nu },\ 0<\nu \leq 1/2,$ consists of shifts.

   In [11] W. Arveson posed a question: to describe the class of cocycle
conjugacy of the flow of shifts on the $W^*$-algebra $\cal M$ defined
by R.T. Powers in [12]. One can see the answer on this question in the
case of $\cal M=B(H)$ in [11]. We investigate cocycle conjugacy
of quasifree automorphisms and endomorphisms semigroups on
the hyperfinite factors of type $II_1$ and $III_{\lambda },\ 0<\lambda <1$.

  In what follows we denote the trace class, the Hilbert Schmidt class,
compact operators and the Hilbert-Schmidt norm by symbols $s_1,s_2,s_{\infty }$
and $||\cdot ||_2$ correspondently.

\medskip

{\bf 2. An extension on $\cal B(H)$ of quasifree automorpisms
of hyperfinite factors $\cal M\subset B(H)$.}

\medskip

Fix a positive operator $R,\ 0<R<I,\ kerR=ker(I-R)=0,$ and
an antiunitary operator $J$ in $K$ and
construct a representation $\pi $ of the algebra $\A(\K \oplus \K)$ in a Hilbert
space $\H =\F (\K)\otimes \F (\K)$ by the formula
$$
\pi (a(f\oplus 0))=a((1-R)^{1/2}f)\otimes \Gamma +1 \otimes a\ (R^{1/2}Jf),
$$
$$
\pi (a(0\oplus f))=a(R^{1/2}f)\otimes \Gamma -1 \otimes a\ ((1-R)^{1/2}Jf),
$$
$f\in \K.$ Here $\Gamma $ is a unitary operator completely defined by
the relations $\Gamma a(f)=-a(f)\Gamma ,\ \Gamma \Omega =\Omega ,\ f\in \K$.
Define the hyperfinite factors $M_R=\pi (\A (\K \oplus 0))''$ and
$M_P=\pi (\A(\K\oplus \K))''=\B(\H)$ and consider two vector states on its,
$$
\omega _R(x)=<\Omega \otimes \Omega ,\pi (x\oplus 0)\Omega \otimes \Omega >,\
x\in \A (\K \oplus 0).
$$
$$
\omega _P(x)=<\Omega \otimes \Omega ,\pi (x)\Omega \otimes \Omega >,\
x\in \A (\K \oplus \K).
$$
Note that
$$
\omega _R(\pi (a^*(f)a(g)))=(f,Rg),\ f,g\in \K,
$$
$$
\omega _P(\pi (a^*(f)a(g)))=(f,Pg),\ f,g\in \K\oplus \K,
$$
where $P=\left (\begin {array} {cc}R&R^{1/2}(I-R)^{1/2}\\
R^{1/2}(I-R)^{1/2}&I-R\end {array}\right )$ is an orthogonal projection
in a Hilbert space $\K\oplus \K$. The state $\omega _R$ is exact
and the state $\omega _P$ is pure
and obtained by the purification procedure (see [19]) from
$\omega _R$. Operators
$$
b(f)=\Gamma \otimes \Gamma \pi (a(0\oplus f)),\
b^*(f)=\pi (a^*(0\oplus f))\Gamma \otimes \Gamma ,\ f\in K,
$$
generate the commutant $M_R'$ and satisfy the relation
$b(f)={\cal J}\pi (a(f\oplus 0)){\cal J},\ f\in K$, where $\cal J$ is
a modular involution on $M_R$ associated with $\omega _R$
(see Appendix).
Let $V$ and $W$ be isometrical operators
in $K$ commuting with $R$. Consider quasifree endomorphisms
$\alpha$ of $M_R$ and $\beta $ of its commutant $M_R'$ obtained by
the lifting of $V$ and $W$: $\alpha (\pi (a(f\oplus 0)))=\pi (a(Vf\oplus 0)),\
\beta (b(f))=b(Wf),\ f\in K$. Consider minimal unitary dilations
$V'$ and $W'$ of the operators $V$ and $W$ acting in the same
Hilbert space $K'$, $K\subset K'$ and a positive contraction
$R'$ in $K'$ such that $V'R=R'V,\ W'R=R'W$ (on the existence of
$R'$ see in [17]).
Determine a quasifree endomorphism
$\theta $ of the $C^*$-algebra generated by $M_R$ and $M_R'$ such that
$\theta |_{M_R}=\alpha ,\ \theta |_{M_R'}=\beta $.

{\bf Theorem 1.} {\it The endomorphism $\theta $
defined by the quasifree lifting of the isometrical operators $V$ and $W$,
can be extended on $M_P$ if and only if the following inclusion
holds, $R^{'1/2}(I-R')^{1/2}(V'-W')\in s_2$. }

{\bf Remark.} {\it The condition of the proposition is sufficient
for the cocycle conjugacy of endomorphic semigroups on $M_R$
obtained by the quasifree lifting of the isometrical
operators $V$ and $W$ included in the semigroups
$\cal V$ and $\cal W$ in a Hilbert space $K$ (see below).}

Proof.

As it was proved by H.Araki (see [9-10]), any quasifree
*-automorphism $\theta '$ given on the $C^*$-algebra
$\pi (A(K'\oplus K'))$ by the formula $\theta '(\pi (a(f\oplus g)))=
\pi (a(V'f\oplus W'g)),\ f,g\in K'$, can be extended on
the factor $M_{P'},\
P'=\left (\begin {array} {cc}R'&R^{'1/2}(I-R')^{1/2}\\
R^{'1/2}(I-R')^{1/2}&I-R'\end {array}\right )$
if and only if
$\left (\begin {array} {cc}V'&0\\
0&W'\end {array}\right )P'-P'
\left (\begin {array} {cc}V'&0\\
0&W'\end {array}\right )\in s_2$. The unitary operator
$\left (\begin {array} {cc}V'&0\\
0&W'\end {array}\right )$
in the space $K'\oplus K'$ satisfies
this condition and correctly define $\theta '$.
The automorphism $\theta '$ has the property $\theta '(\Gamma \otimes \Gamma )=
\Gamma \otimes \Gamma $, such that $\theta '(b(f))=b(W'f),\ f\in K'$.
The factor $M_P\subset M_{P'}$ is invariant under
the action of $\theta '$. Considering the restriction we obtain
$\theta '|_{M_P}=\theta $.
Note that in the case of $V=W$, the endomorphism $\theta $ is
a regular extension of $\alpha $ in the sense of [8]. $\triangle $

Now let $(V_t)_{t\in R_+}$ be a $C_0$-semigroup of isometrical operators
in $K$. Then a family of quasifree endomorphisms on
$M_R$ defined by the formula $\alpha _t(\pi (a(f\oplus 0)))=
\pi (a(V_tf\oplus 0)),\ f\in K,\ t\in R_+,$ is a $E_0$-semigroup.
Involve an expanding family of Hilbert spaces $({\cal H}_t)_{t\in R_+}$
embedded in $F(K)$ such that
${\cal H}_t$ is generated by all vectors
$a^{\# }(f_1)a^{\# }(f_2)\dots a^{\# }(f_n))\Omega ,\ f_i\in kerV_t^*,\
1\leq i\leq n$, where $a^{\# }=a^*$ or $a$. The family $({\cal H}_t)_{t\in R_+}$
is a product-system of Hilbert spaces (see [8]).
Consider a product-system $K_t={\cal H}_t\otimes J{\cal H}_t,\ t\in R_+$.
Let $(\theta _t)_{t\in R_+}$ be a regular extension $(\alpha _t)
_{t\in R_+}$, then
$$
K_t=\theta _t(|\Omega \otimes \Omega ><
\Omega \otimes \Omega |){\cal H},\ t\geq 0.
$$
By this way, the product-system associated
with the regular extension is obtained by the doubling
of the product-system associated with the initial semigroup.
Note that in the common case (see [22])
$\theta _t(P)$, where $P$ is a one-dimensional projection,
is not obliged to be a
monotonous increasing family of projections.
Such situation can characterize the
complete compatability with the exact state
( in this case it is $\omega _R$).

\medskip

{\bf 3. Inner *-automorphisms and the cocycle conjugacy on the hyperfinite
factor $\cal M\it _{R}$.}

\medskip

   In [9] it was proved that a quasifree derivation $\delta $
of $\A (\K)$ acting on the generating elements by the
formula $\delta (a(f))=a(df),\ f\in dom\ d$, where $d$ is some scewhermitian
operator, is inner iff $d\in s_1$. Then the automorphism
obtained by the quasifree lifting of an unitary $e^d,d\in s_1$, is
inner. Notice that $e^d-I\in s_1$. On the other side every
inner automorphism $\alpha $ can not be represented in the form
$\alpha =e^{\delta }$, where $\delta $ is some inner derivation. In the
following theorem we give the necessary and sufficient condition
of an innerness of $B_{R}(W)$ in the terms of $W,\ WR=RW$.

   {\bf Theorem 2.} {\it The quasifree automorphism $B_{R}(W)$ of
the hyperfinite factor $\cal M\it _{R}$ is inner iff
$R^{1/2}(I-R)^{1/2}(W-I)\in s_2$.}

   {\bf Remark.} {\it The results of [18] yelds a sufficiency and
a necessity of the condition of theorem 2 in the case of a pure point
spectrum of $W$. Thus we need to prove a necessity of one.}

   Proof of theorem 2 (necessity).

   The quasifree automorphism $B_{R}(W)$ has
a form $B_{R}(W)(\cdot )=\cal U\it \cdot \cal U\it ^*,\
\cal U\in M\it _{R}$ by the condition. In the appendix we show that
the generating elements of commutant ${\cal M}_{R}'$ are
$\bf 1,\it b(f)=\Gamma \otimes \Gamma \pi (a(0\oplus f)),\
b^*(f)=\pi (a^*(0\oplus f))\Gamma \otimes \Gamma ,\ f\in \K$.
Thus ${\cal U}\it \pi (a(0\oplus f)){\cal U}\it ^*={\cal U}
\Gamma \otimes \Gamma {\cal U}^*\Gamma \otimes \Gamma \pi
(a(0\oplus f)),\ f\in \K$. Let us show that
${\cal U}\Gamma \otimes \Gamma {\cal U}^*\Gamma \otimes \Gamma =1$.
Note that $(\Gamma \otimes \Gamma )^2=I,\ (\Gamma \otimes \Gamma )^*=
\Gamma \otimes \Gamma $ therefore ${\cal U}\Gamma \otimes \Gamma
{\cal U}^* \Gamma \otimes \Gamma $ equals $\bf 1$ or $-\bf 1$.
By this way
${\cal U}\pi (a(0\oplus f)){\cal U}^*=\ \pi (a(0\oplus f))$ or
$-\pi (a(0\oplus f)),\ f\in \K$. In the first case the quasifree automorphism
$B(W\oplus I)$ is unitary implementable, in the second case the quasifree
automorphism $B(W\oplus (-I))$ is that. Therefore by theorem 1,
$R^{1/2}(I-R)^{1/2}(W-I)\in s_2$ or
$R^{1/2}(I-R)^{1/2}(W+I)\in s_2$. In the first case the theorem is proved.
Suppose $R^{1/2}(I-R)^{1/2}(W+I)\in s_2$.
Then $R^{1/2}(I-R)^{1/2}(-W-I)\in s_2$ and the automorpism
$B_{R}(-W)$
is inner. Thus the automorphism $B_{R}(-I)=B_{R}(-W)B_{R}(W)$ is
inner. It is a contradiction by [18]. Thus we proved
${\cal U}\pi (a(0\oplus f)){\cal U}^*=\pi (a(0\oplus f))$.
Therefore the automorphism $B(W\oplus I)$ of the $C^*$-algebra
$\pi (\A (\K \oplus \K))$ obtained by the lifting of
$W\oplus I$ is unitary implementable and
$R^{1/2}(I-R)^{1/2}(W-I)\in s_2$ by theorem 1. $\triangle $

   Let $(U_t)_{t\geq 0}$ and $(V_t)_{t\geq 0}$ be $C_0$-semigroups
of unitaries in a Hilbert space $\K$ commuting with
the operator $R$.

{\bf Theorem 3.} {\it The quasifree semigroups $(B_{R}(U_t))_{t\geq 0}$
and $(B_{R}(V_t))_{t\geq 0}$ on the hyperfinite factor ${\cal M}_{R}$
are cocycle conjugate iff $R^{1/2}(I-R)^{1/2}(U_t-V_t)\in s_2,\ t\geq 0$.}

   Proof of theorem 3 is based on theorem 1, theorem 2 and one result
of [20] that we formulate in the following lemma:

{\bf Lemma.} {\it Let $\sigma $-week continuous groups of *-automorphisms
$\alpha =(\alpha _t)_{t\in\bf R}$ and $\beta =(\beta _t)_{t\in \bf R}$ on
the $W^*$-factor $\cal M$ having separable predual $\cal M\it _*$ be such that
*-automorphisms $\beta _{-t}\alpha _t$ are inner for all $t\in \bf R$.
Then the group $\alpha $ and $\beta $ are cocycle conjugate.}

   Proof of theorem 3.

   Necessity. Let the semigroups $(B_{R}(U_t))_{t\geq 0}$ and
$(B_{R}(V_t))_{t\geq 0}$ be cocycle conjugate. Then
$B_{R}(U_t)(\cdot )=W_tB_{R}(V_t)(\cdot )W_t^*,\ W_t\in {\cal M}_{R},\
t\geq 0$. Fix $t\geq 0$. The automorphism $B(V_t\oplus V_t)$ is unitary
implementable
by theorem 1. Therefore the automorphism $B(U_t\oplus V_t)(\cdot )=
W_tB(V_t\oplus V_t)(\cdot )W_t^*$ is unitary implementable too. The result follows
from theorem 1.

   Sufficiency. The result follows from theorem 2 and the lemma. $\triangle $

\medskip

{\bf 4. The cocycle conjugacy of quasifree endomorphisms semigroups.}

\medskip

   Let $U=(U_t)_{t\geq 0}$ and $V=(V_t)_{t\geq 0}$ be $C_0$-semigroups of isometries
in a Hilbert space $\K$ commuting with the operator $R$.
Then there are $C_0$-semigroups of unitaries $U'=(U_t')_{t\geq 0}$
and $V'=(V_t')_{t\geq 0}$ in a Hilbert space $\K',\ \K\subset \K',$
being minimal unitary
dilations of semigroups $U$ and $V$ and the positive contraction $R'$
commuting with $U',V'$ and satisfying the relation
$U'_tR=R'U_t,\ V'_tR=R'V_t,\ t\geq 0$.

   {\bf Definition.} {\it Two $C_0$-semigroups of isometries $U$ and $V$ are called
to be approximating each other if $U_t'-V_t'\in s_2$ and
$U_t'V_t^{'*}|_{\K'\ominus \K}=I,\ t\geq 0$.}

   {\bf Theorem 4.} {\it Let a $C_0$-semigroup of isometries $U$ approximates
a $C_0$-semigroup of isometries $V$. Then the quasifree semigroup
$(B_{R}(U_t))_{t\geq 0}$
is cocycle conjugate to the quasifree semigroup $(B_{R}(V_t))_{t\geq 0}$.}

   Proof of theorem 4.

   Theorem 3 leads to an existence of a cocycle $(\cal W\it _t)_{t\geq 0}$
such that $B_{R'}(U_t')(\cdot )=\cal W\it _tB_{R'}(V_t')(\cdot )\cal W\it _t^*,\
\cal W\it _t\in {\cal M}_{R'}(\K '),\ t\geq 0$. We show that
the condition $U_t'V_t^{'*}|_{\K'\ominus \K}=I$ implies $\cal W\it _t\in
\cal M\it _{R}(\K)$. Note that if this condition holds
a family of unitaries $W_t=U_t'V_t^{'*}|_{\K},\ W_t-I\in s_2,\ t\geq 0,$ is
correctly defined. The quasifree lifting of $(W_t)_{t\geq 0}$ determines
a family of inner automorphisms $B_{R}(W_t)(x)=\cal W\it '_tx
\cal W\it _t^{'*},\ \cal W\it _t',x\in \cal M\it _{R}(\K),\ t\geq 0$.
We constructed $({\cal W}_t')_{t\geq 0}$ such that
$\cal W\it _t'x\cal W\it _t^{'*}=\cal W\it _tx\cal W\it _t^*,\
x\in \cal M\it _{R'}(\K'),\ t\geq 0$. Therefore $\cal W\it _t=
e^{ic(t)}\cal W\it _t' ,\ c(t)\in {\bf R},\ t\geq 0,$ and
$\cal W\it _t\in \cal M\it _{R}(\K)$.

\medskip

{\bf 5. Continuous semigroups of isometries in a Hilbert space.}

\medskip

It is useful to remind that an isometry $V$ in a Hilbert space $\K$ is called
completely nonunitary if there is no subspace $\K _0\subset \K$ reducing $V$
to an unitary. Every $C_0$-semigroup of completely nonunitary
isometries $(V_t)_{t\geq 0}$ is unitary equivalent to its model that is
a $C_0$-semigroup of shifts $(S_t)_{t\geq 0}$ acting in a Hilbert space
$\K '=H\otimes L_2(0,+\infty )$ by the formula $(S_tf)(x)=f(x-t)$ for
$x>t,\ (S_tf)(x)=0$ for $0<x<t,\ f\in \K '$.
Here $H$ is some Hilbert space of the dimension $n$ equal to the deficiency
index of the generator of $(V_t)_{t\geq 0}$ (we call $n$ by the deficiency
index of $(V_t)_{t\geq 0}$ in the following).
Let $(V_t)_{t\geq 0}$ be a $C_0$-semigroup of isometries in a Hilbert space $\K$
and $\K=\K _0\oplus \K _1$ be the Wold decomposition of $\K$, where $\K _0$ reduces
$(V_t)_{t\geq 0}$ to a $C_0$-semigroup of unitaries and $\K _1$ reduces
$(V_t)_{t\geq 0}$
to a $C_0$-semigroup of completely nonunitary isometries. We shall say
$(V_t)_{t\geq 0}$ satisfies the condition $N$ if the $C_0$-semigroup of
unitaries $(V_t|_{\K _0})_{t\geq 0}$ is uniformly continuous that is
$||V_t|_{\K _0}-V_s|_{\K _0}||\to 0,$ $t\to s,\ s,t\geq 0$.

{\bf Theorem 5 ([1-3,5-6]).} {\it Let a $C_0$-semigroup of isometries $(V_t)_{t\geq 0}$
in a Hilbert space $\K$ with a deficiency index $n>0$ satisfy
the condition $N$.

   Then there is a $C_0$-semigroup of completely nonunitary isometries
$(S_t)_{t\geq 0}$ with a deficiency index $n$ approximating $(V_t)_{t\geq 0}$
in the sence of the definition of part 4.}

   In the proof of the theorem we use the complex analysis in
the Hardy space (see [21]). Let $(\lambda _k)_{1\leq k\leq N},\
N\leq +\infty$, be a system of complex numbers satisfying the following
properties,
$$
Re\lambda _k<0,\ |Im\lambda _k|<R,\ 1\leq k\leq N,\
\sum \limits _{k=1}^{N}|Re\lambda _k|<+\infty ,
\eqno (1)
$$
where $R$ is some positive number. In this case the formula
$B(\lambda )=\prod \limits _{k=1}^{N}\frac {\lambda +\overline \lambda _k}
{\lambda -\lambda _k},\ \lambda \in \bf C$, defines an analitic function
being regular in the semiplane $Re\lambda >0$ and equaling the unit by module
on the imaginary axis. The function $B(\lambda )$ is called the Blaschke product.

   {\bf Proposition 1.} {\it Let the numbers $(\lambda _k)_{1\leq k\leq N}$
satisfy the condition (1). Then the Blaschke product $B(\lambda )$
constructed of $(\lambda _k)_{1\leq k\leq N}$ can be estimated as
follows:
$$
|B(\lambda )|<C_1,\ |\lambda |>C_2,\ B(\lambda )=1-\frac {C_3}{\lambda }+
o(\frac {1}{\lambda }),\ |\lambda |\to +\infty ,\ \lambda \in \bf C,
$$
where $C_1,C_2$ and $C_3$ are some positive constants.}

   Proof.

$$
ln B(\lambda )=\sum \limits _{k=1}^{N}ln\frac {1+\frac {\overline {\lambda _k}}
{\lambda }}{1-\frac {\lambda _k}{\lambda }}=-\frac {2s}{\lambda }+
o(\frac {1}{\lambda }),\ |\lambda |\to +\infty ,
$$
where $s=-\sum \limits _{k=1}^{N}Re\lambda _k$, $0<s<+\infty $ by the condition
(1). $\triangle $

   Let us consider a Hilbert space $\K=L_2(0,+\infty )$. Let $P_{[t_1,t_2]}$
designate a projection on a subspace of $\K$ consisting of functions
$f(x)=0$ for $0<x<t_1,\ t_2<x<+\infty $. Let us define
an isometry $\Theta $ acting in $\K$ by the formula $\Theta =
\cal F\it ^{-1}B\cal F$, where $\cal F$ and $B$ are the Fourier transformation
and an operator of the multiplication by the Blaschke product correspondently.

   {\bf Proposition 2.} {\it Let the conditions of proposition 1 be hold.
$\\ $Then $\Delta _{t,\delta }=P_{[t,t+\delta ]}\Theta P_{[t,t+\delta ]}-
P_{[t,t+\delta ]}\in s_2,\ 0<t,\delta <+\infty ,\ ||\Delta _{t,\delta }||_2=
O(\delta ^{1/2}),\ \delta \to 0.$}

   Proof.

   Fix $t\geq 0$. Let $\mu _{k,\delta }=-\frac {1}{2|k|}+i\frac {2\pi k}
{\delta },\ k\in {\bf Z},\ \delta >0$.
Let us consider a family of functions $f_{k,\delta }(x)=
\frac {(-2Re\mu _{k,\delta })^{1/2}}{(e^{2Re\mu _{k,\delta }t}-
e^{2Re\mu _{k,\delta }(t+\delta )})^{1/2}}e^{\mu _{k,\delta }x},\
t<x<t+\delta ,\ f_{k,\delta }(x)=0,
\ 0<x<t,\ t+\delta <x<+\infty ,\ k\in \bf Z$. The family
$(f_{k,\delta })_{k\in \bf Z}$ is the Riesz basis of a Hilbert space
$H=P_{[t,t+\delta ]}\K$ (see [21]). So there is a bounded operator
$V$ having bounded revers in $\K$ such that the family $(Vf_{k,\delta })_{k\in \bf Z}$
is an orthogonal basis of $H$. Therefore to prove proposition 2 it is sufficient
to prove a convergence of the series $\sum \limits _{k\in \bf Z}
||\Delta _{t,\delta }f_{k,\delta }||^2$ and to investigate its dependence on
$\delta $. Let $f_{k,\delta }=
\frac {(-2Re\mu _{k,\delta })^{1/2}}{e^{2Re\mu _{k,\delta }t}-
e^{2Re\mu _{k,\delta }(t+\delta ))^{1/2}}}e^{\mu _{k,\delta }x},\ f^{(2)}_{k,\delta }=P_{[t+\delta ,+\infty ]}f^{(1)}
_{k,\delta },\ t<x<+\infty ,\ f^{(1)}_{k,\delta }=f^{(2)}_{k,\delta }=0,\
0<x<t,\ k\in \bf Z$. Then $f_{k,\delta }=f^{(1)}_{k,\delta }-f^{(2)}_{k,\delta },
\ k\in \bf Z$ and $||\Delta _{t,\delta }f_{k,\delta }||^2\leq
2(||\Delta _{t,\delta }f^{(1)}_{k,\delta }||^2+||\Delta _{t,\delta }
f^{(2)}_{k,\delta }||^2),\ k\in \bf Z$. By this way,
$$
||\Delta _{t,\delta }f^{(i)}_{k,\delta }||^2=
2(||f^{(i)}_{k,\delta }-Re(\Theta f^{(i)}_{k,\delta },f^{(i)}_{k,\delta })),\
k\in {\bf Z},\ i=1,2,
$$
$$
||f_{k,\delta }^{(1)}||^2=\frac {e^{Re\mu _{k,\delta }t}}{(e^{Re\mu _{k,\delta }t}-
e^{Re\mu _{k,\delta }(t+\delta )})}=-\frac {1}{2Re\mu _{k,\delta }\delta }+o(1),
\ |k|\to +\infty ,
$$
$$
||f_{k,\delta }^{(2)}||^2=\frac {e^{Re\mu _{k,\delta }(t+\delta )}}
{(e^{Re\mu _{k,\delta }t}-e^{Re\mu _{k,\delta }(t+\delta )})}=
-\frac {1}{2Re\mu _{k,\delta }\delta }+o(1),\ |k|\to +\infty .
$$
Using the Laplace transformation technics one can obtain
$$
(\Theta f_{k,\delta }^{(i)},f_{k,\delta }^{(i)})=B(\mu _{k,\delta })
||f_{k,\delta }^{(i)}||^2,\ i=\overline {1,2}.
$$
It follows from proposition 1 that
$$
||\Delta _{t,\delta }f_{k,\delta }^{(i)}||^2=
||f_{k,\delta }^{(i)}||^2(1-ReB(\mu _{k,\delta }))=\frac {1}
{2\delta |\mu _{k,\delta }|^2}+o(\frac {1}{2\delta |\mu _{k,\delta }|^2})=
$$
$$
\frac {\delta }{8\pi ^2k^2}+o(\frac {\delta }{k^2}),\
|k|\to +\infty ,\ \delta \to 0,\ i=\overline {1,2},
$$
and
$$
||\Delta _{t,\delta }||^2_2\leq C\sum \limits _{k\in \bf Z}
||\Delta _{t,\delta }f_{k,\delta }||^2\leq
2C\sum \limits _{k\in \bf Z}(||\Delta _{t,\delta }f_{k,\delta }^{(1)}||^2+
||\Delta _{t,\delta }f_{k,\delta }^{(2)}||^2)=O(\delta ),\ \delta \to 0,
$$
where $C$ is some positive constant, that implies
$
||\Delta _{t,\delta }||_2=O({\delta }^{1/2}),\
\delta \to 0.\ \triangle
$

   Proof of theorem 5.

   Let the operator $d$ be a generator of uniformly continuous semigroup
of unitaries $(U_t)_{t\geq 0}$ being unitary part of some
semigroup of isometries. In accordance with the Von Neumann theorem, given
a scewhermitian operator $d$ there is a bounded scewhermitian operator
$D\in \sigma _2$ such that the scewhermitian operator $d+D$ has a purely point
spectrum. So the operator $d+D$ is a generator of some uniformly
continuous semigroup of unitaries $(V_t)_{t\geq 0}$ having a purely point spectrum.
The semigroups $(U_t)_{t\geq 0}$ and $(V_t)_{t\geq 0}$
are known to be bonded by the relation $V_t-U_t=\int \limits _{0}^{t}
U_{t-s}DV_sds,\ t\geq 0$, therefore $V_t-U_t\in \sigma _2,\ t\geq 0$ and
$||V_{t+\delta }-U_{t+\delta }-V_t+U_t||_2=O(\delta ),\ \delta \to 0,\
t\geq 0$.

   For the semigroup $(V_t)_{t\geq 0}$ is uniformly continuous its spectrum
lies in a circle of radius $R>0$ in the complex plane. Let $(i\mu _k)_{1\leq k\leq N}$
be eigenvalues of the generator of $(V_t)_{t\geq 0}$ numbering in order of
decreasing of its modules. By this way, $|\mu _k|<R,\ 1\leq k\leq N\leq +\infty $.
Let complex numbers $(\lambda _k)_{1\leq k\leq N}$ be such that
$Im\lambda _k=-\mu _k,\ 1\leq k\leq N,$ and the real parts collected in
the accordance of
the condition (1). The condition (1) allows to define the Blaschke product
associated with $(\lambda _k)_{1\leq k\leq N}$ (see above).

   In what follows we suppose the deficiency index to be $1$.
It will prove the theorem because every $C_0$-semigroup of isometries
having nonzero deficiency index decomposes in
an orthogonal sum of a $C_0$-semigroup of isometries having a deficiency
index $1$ and, probably, a $C_0$-semigroup of completely
nonunitary isometries. We shall prove the existence of a $C_0$-semigroup of
isometries in a Hilbert space $\K=L_2(0,+\infty )$ that is unitary equivalent
to given a semigroup of isometries with a deficiency index $1$
and a purely point spectrum of its unitary part consisting of numbers
$(i\mu _k)_{1\leq k\leq N}$.

   Let $(S_t)_{t\geq 0}$ be a  $C_0$-semigroup of shifts in a Hilbert
space $\K=L_2(0,+\infty )$. Let us consider a family of functions
$f_n(x)=(-2Re\lambda _n)^{1/2}e^{\lambda _nx},\ 1\leq n\leq N$.
The condition (1) implies an uncomleteness of the system $(f_n)_{1\leq n\leq N}$
in $\K$ ( the condition of the convergence of the Blaschke product ).
Thus a subspace $\K _1$ being a linear envelope of $(f_n)_{1\leq n\leq N}$
does not coincide with $\K$ and defines a subspace $\K _0$ being invariant
under an action of $(S_t)_{t\geq 0}$ and completely describing by the condition
of an orthogonality to all functions $f_n$ such that $\K=\K _0\oplus \K _1$
and the isometry $\Theta :\K\to \K,\ \Theta =\cal F\it ^{-1}B\cal F,\it \
\Theta S_t=S_t\Theta ,\ t\geq 0$, where $\cal F$ and $B$ are the Fourier
transformation
and an operator multiplying by the Blaschke product, defines $\K _0$ and
$\K _1$ by the formula $\K _0=\Theta \K$. The semigroup $(S_t)_{t\geq 0}$ is
intertwined by the operator $\Theta $
with its restriction on the subspace $\K _0$: $S_t|_{\K _0}\Theta =\Theta S_t,\
t\geq 0$. The isometric operator $\Theta :\ \K\to \K$ sets an
unitary map $\K\to Ran\Theta =\K _0$. Hence the semigroups $(S_t)_{t\geq 0}$ and
$(S_t|_{h_0})_{t\geq 0}$ are unitarily equivalent such that the deficiency
index of the semigroup $(S_t|_{\K _0})_{t\geq 0}$ consisting
of completely nonunitary isometries in $\K _0$ coincides with the deficiency
index of the semigroup $(S_t)_{t\geq 0}$ and equals $1$.

   Let the system of functions $(g_n)_{n\in \bf N}$ be obtained by a succesive
orthogonalization of the system $(f_n)_{n\in \bf N}$. Let us determine
a $C_0$-semigroup $(V_t)_{t\geq 0}$ of isometries in $\K$ as follows
$$
V_t|_{\K _0}=S_t|_{\K _0},\ V_tg_n=e^{iIm\lambda _nt}g_n,\ t\geq 0,n\in \bf N.
\eqno (2)
$$

   We shall show that for isometries $V_t,\ t\geq 0,$ describing in (2)
the following conditions hold, $V_tS_t^*-P_{[t,+\infty ]}\in \sigma _2,
||V_t-S_t||_2\leq ||V_tS_t^*-P_{[t,+\infty )}||_2=O(t^{1/2}),\ t\to 0$.

   Fix $t>0$. We need to prove a convergence of the series
$\sum \limits _{n=1}^{+\infty }||(V_tS_t^*-P_{[t,+\infty )})f_n||^2$
for some orthogonal basis $(f_n)_{n\in \bf N}$ of $\K$.
Choose for this purpose an arbitrary addition of the system $(S_tg_n)
_{n\in \bf N}$ up to an orthogonal basis of $\K$.

   Notice that $V_tS_t^*-P_{[t,+\infty )}=P_{[0,t]}V_tS_t^*+
(P_{[t,+\infty )}V_tS_t^*-P_{[t,+\infty )})$ and
$(V_tS_t^*-P_{[t,+\infty )})|_{\K _t}=0,$ where $\K _t$ designates an orthogonal
addition of a linear envelope of vectors $(S_tg_n)_{n\in \bf N}$.
An element $S_t^*g_n$ belongs to a linear envelope of elements $g_i,\ i=1,n$,
such that $(S_t^*g_n,g_n)=e^{\lambda _nt},\ 1\leq n\leq N$.
By this way,

$$
||(P_{[t,+\infty )}V_tS_t^*-P_{[t,+\infty )})S_tg_n||^2=1+
||P_{[t,+\infty )}g_n||^2-2Re(V_tg_n,S_tg_n)<
$$
$$
2(1-Re(V_tg_n,S_tg_n))=2(1-e^{-Re\lambda _n}),\ n\in \bf N.
\eqno (3)
$$

   For the operator $P_{[0,t]}V_tS_t^*$ we have the estimate using
the Bessel inequality:

$$
||P_{[0,t]}V_tS_t^*S_tg_n||^2=||P_{[0,t]}g_n||^2=1-
(P_{[t,+\infty )}g_n,g_n)<1-|(S_t^*g_n,g_n)|^2
$$
$$
\ \ \ \ \ \ \ \ \ \ \ \ \ \ \ \ \ \ \ \ \ \ \ \ \ \ \ \ \ \ \ \ \ \ \ \ \ \
\ \ \ \ \ \ \ \ \ \ \ \ \ \ \ \ \ \ \ =1-e^{-2Re\lambda _nt},\ n\in \bf N.
\eqno (4)
$$

   It follows from $(3),(4)$ and $(1)$ that
$\\ \sum \limits _{i=1}^{+\infty }||(V_tS_t^*-P_{[t,+\infty )})f_i||^2=
\sum \limits _{n=1}^{+\infty }(||(P_{[t,+\infty )}V_tS_t^*-P_{[t,+\infty )})
S_tg_n||^2+||(P_{[0,t]}V_tS_t^*S_tg_n||^2)=\\ $
$\sum \limits _{n=1}^{+\infty }(-4Re\lambda _nt+o(Re\lambda _n))=O(t),\
t\to 0$. Hence
$||V_tS_t^*-P_{[t,+\infty )}||_2=O(t^{1/2}),\ t\to 0$.

   Notice that $V_t-S_t=V_t|_{\K _1}-P_{\K _1}S_t|_{\K _1},\ t\geq 0,$ and
the $C_0$-semigroup $(V_t|_{\K _1})_{t\geq 0}$ is uniformly continuous by the
condition.
The condition $V_t-S_t\in s_2,\ t\geq 0,\ ||V_t-S_t||_2=O(t^{1/2}),\
t\to 0,$ and the uniform continuity of $(V_t|_{\K _1})_{t\geq 0}$ implies that
the family $(V_t-S_t)_{t\geq 0}$ is continuous in $||\cdot ||_2$. In fact,
$||(V_{t+\delta }-S_{t+\delta }-V_t+S_t||_2=||P_{\K _1}(V_{t+\delta }-
S_{t+\delta }-V_t+S_t)P_{\K _1}||_2\leq ||P_{\K _1}(V_{\delta }-I)(V_t-S_t)
P_{\K _1}||_2
+||(V_{\delta }-S_{\delta })V_t||_2+||(V_t-S_t)(V_{\delta }-S_{\delta })||_2\leq
||P_{\K _1}(V_{\delta }-I)P_{\K _1}||||V_t-S_t||_2+||V_{\delta }-S_{\delta }||_2
(1+||V_t-S_t||_2)\to 0,\ \delta \to 0.$

   Now we get for operators $\Delta _t=V_t-S_t,\ t\geq 0$ the following
estimations:
$$
\Delta _t\in s_2,\ ||\Delta _{t+\delta }-\Delta _t||_2\to 0,\
||\Delta _{\delta }||_2=O(\delta ^{1/2}),\ \delta \to 0,\ t\geq 0.
\eqno (5)
$$

   To complete the proof we need to show that there are $C_0$-semigroups
$(V_t')_{t\geq 0}$ and $(S_t')_{t\geq 0}$ being unitary dilations of
$(V_t)_{t\geq 0}$ and $(S_t)_{t\geq 0}$ which satisfy the conditions
$V_t'-S_t'\in s_2,\ V_t'S_t^{'*}|_{\K '\ominus \K}=I,\ t\geq 0$.
Let us define unitary dilations in a Hilbert space $\K '=\K\oplus \K$ by
the formula
$$
S_t'(f\oplus g)(x)=((S_tf)(x)+(P_{[0,t]}g)(t-x))\oplus (S_t^*g)(x),
$$
$$
V_t'(f\oplus g)(x)=((V_tf)(x)+(\Theta (P_{[0,t]}g)(t-\cdot ))(x))
\oplus (S_t^*g)(x),
$$
$$
x,t\geq 0,\ f,g\in \K.
$$
Now the result follows from (5) and proposition 2. $\triangle $

\medskip

{\bf 6. The class of the cocycle conjugacy of the quasifree flow of shifts
on the hyperfinite factor ${\cal M}_{\nu }$.}

\medskip

   Let $V=(V_t)_{t\geq 0}$ be a $C_0$-semigroup of isometries having
uniformly continuous unitary part (see part 5). Then theorem 5 implies
the existence of a $C_0$-semigroup of completely nonunitary isometries
$S=(S_t)_{t\geq 0}$ in $\K$ approximating $V$ in the sence of part 4.
Semigroups $V$ and $S$ have same deficiency indeces.
If $(S_t)_{t\geq 0}$ is a $C_0$-semigroup of completely nonunitary isometries,
then the quasifree semigroup $(B_{\nu }(S_t))_{t\geq 0}$ consists of shifts.
It follows from theorem 4 that the following assertion holds.

{\bf Theorem 6 ([2-5]).} {\it Let a $C_0$-semigroup of isometries
$(V_t)_{t\geq 0}$ with the deficiency index $n>0$ have uniformly continuous
unitary part.

   Then there is a $C_0$-semigroup of completely nonunitary isometries
$(S_t)_{t\geq 0}$ with the deficiency index $n$ such that
the quasifree semigroup $(B_{\nu }(V_t))_{t\geq 0}$ is cocycle conjugate to
the flow of shifts $(B_{\nu }(S_t))_{t\geq 0}$.}

Theorem 6 describes the class of the cocycle conjugacy
of the quasifree flow of shifts on the hyperfinite factor $\cal M\it _{\nu }$
so one can reformulate it as follows,

{\bf Theorem 6'.} {\it Let $V=(V_t)_{t\geq }$ be a $C_0$-semigroup of isometries
in a Hilbert space $\K$ having uniformly continuous unitary part
and $S=(S_t)_{t\geq 0}$ be a $C_0$-semigroup of completely nonunitary isometries
in $\K$ with the same deficiency index as in $V$.

   Then there is an unitary $U$ in $\K$ such that the quasifree semigroup
$(B_{\nu }(UV_tU^*))_{t\geq 0}$ is cocycle conjugate to
the flow of shifts $(B_{\nu }(S_t))_{t\geq 0}$.}

Notice that there is an analogue of the result of theorem 6 for discrete quasifree
semigroups, see [5,7].

\medskip

{\bf Acknowledgements.} {\it The author would like to thank
A.V. Bulinskij for the initiation of this work and useful
discussions.}

\medskip

{\bf Appendix. The commutant of
${\cal M}_{R}$.}

\medskip

Note that the operators
$$
b(f)=\Gamma \otimes \Gamma \pi (a(0\oplus f)),
b^*(f)=\Gamma \otimes \Gamma \pi (a^*(0\oplus f)),\ f\in \K,
$$
belong to a commutant of the hyperfinite factor $\cal M\it _{R}$. Therefore
its generate the $W^*$-algebra $\cal N\it \subset \cal M\it _{R}'$.
The formula $Sx\Omega \otimes \Omega =x^*\Omega \otimes \Omega ,
\ x\in \cal M\it _{\nu }$ correctly defines an antilinear operator $S$ in
$\cal H$ (see [4]).
Let $S=\cal J\it \Delta ^{1/2}$ be a polar decomposition of $S$. Then
an antiisometrical part $\cal J$ of the operator $S$ is antiunitary and
is called a modular involution of $\cal M\it _{R}$.
The linear (unbounded) positive operator $\Delta $ is called a modular
operator (see [14]). Simple calculation gives the following formula
for $\cal J$,
$$
\cal J\it f_1\Lambda ...\Lambda f_n\otimes g_1\Lambda ...\Lambda g_m=
Jg_m\Lambda ...\Lambda Jg_1\otimes Jf_n\Lambda ...\Lambda Jf_1,\
\cal J\it \Omega \otimes \Omega =\Omega \otimes \Omega .
$$
By this way,
$$
\cal J\it \pi (a(f\oplus 0))\cal J=\it b^*(f),\ \ f\in \K,
$$
therefore $\cal JM_{R}J=M_{R}'$ and ${\cal N=M}_{R}'$.

\medskip

{\bf References}

[1]. G.G. Amosov, {\it On approximation of continuous semigroups of isometries},
in Algebra and analysis. Proceed. of Conf. dedicated to B.M. Gagaev,
Kazan, 1997, P. 17-18.

[2]. G.G. Amosov, {\it On cocycle conjugacy classes of quasifree K-systems}, in
Some probl. of fund. and appl. math., Moscow Inst. Phys. and Tech, 1997,
P. 4-16.

[3]. G.G. Amosov, {\it Cocycle conjugacy classes of Powers' semiflows of
shifts}, in Internat. Congr. Math. Abst. of Short Commun. and Post.,
Berlin, 1998, P. 113-114.

[4]. G.G. Amosov, {\it Index and cocycle conjugacy of the endomorphism semigroups
on $W^*$- factors},
Ph.D. Thesis, Moscow Inst. Phys. and Tech., 1998.

[5] G.G. Amosov, {\it Cocycle perturbation of quasifree algebraic
K-flow leads to required asymptotic dynamics of associated
completely positive semigroup}, Infin. Dimen. Anal., Quantum. Prob.
and Rel. Top. 3 (2000) 237-246.

[6]. G.G. Amosov, {\it On approximation of isometrical semigroups in
a Hilbert space}, Izvest. Vysch. Uchebn. Zaved. Mat.
(J. Russ. Math.) No 2 (2000) 7-12.

[7]. G.G. Amosov, {\it Approximation by modulo $s_2$ of isometrical
operators and cocycle conjugacy of endomorphisms on the CAR algebra},
Fund. Appl. Math. (to appear)

[8]. G.G. Amosov, A.V. Bulinskij, {\it The Powers-Arveson index for dynamical
semigroups on $W^*$-algebras}, Mathematical Notes 62 (1997) 781-783.

[9]. H. Araki, {\it Bogoliubov automorphisms and Fock representations of
canonical anticommutation relations}, Contemp. Math. 62 (1985) 21-141.

[10]. H. Araki, {\it On quasifree states of CAR and Bogoliubov automorphisms},
Publ. RIMS. 6 (1971) 385-442.

[11]. W. Arveson, {\it Continuous analogues of Fock space},
Mem. AMS 409 (1989) 1-66.

[12]. R.T. Powers, {\it An index theory for semigroups of *-endomorphisms of
$\cal B(H)$ and $\rm II_1$ factor}, Canad.J.Math. 40 (1988) 86-114.

[13]. A.V. Bulinskij, {\it Algebraic K-systems and Powers's semiflows of shifts},
Uspekhi Mat. Nauk 51 (1996) 145-146.

[14]. U. Bratteli, D. Robinson, Operator algebras and quantum statistical
mechanics I, Springer-Verlag, 1982.

[15]. H. Araki, {\it Expansional in Banach algebras}, Ann. Sci. Ecole Norm. Sup.
6 (1973) 67-84.

[16]. E.B. Davies, {\it Irreversible dynamics of infinite fermion systems},
Commun. Math. Phys. 55 (1977) 231-258.

[17]. D. Evans, {\it Completely positive quasifree maps on the CAR algebra},
Commun. Math. Phys. 70 (1979) 53-68.

[18]. T. Murakami, S. Yamagami, {\it On types of quasifree representations of
Clifford algebras}, Publ. RIMS 31 (1995) 33-44.

[19]. R.T. Powers, E. Stormer, {\it Free states of canonical anticommutation
relations}, Commun. Math. Phys. 16 (1970) 1-33.

[20]. A. Connes, {\it Une classification des facteurs de type III}, Ann. Sci.
Ecole Norm. Sup. 6 (1973) 133-252.

[21]. N.K. Nikolski, Treatise on the shift operator, Springer-Verlag, 1986.

[22]. B.V.R. Bhat, {\it An index theory for quantum dynamical semigroups},
Trans. AMS 348 (1996) 561-583.

\end {document}